\documentclass[10pt,letterpaper,conference]{IEEEconf}  
 \usepackage{etoolbox}

\makeatletter
\patchcmd{\@begintheorem}{\textit}{\textbf}{}{}
\makeatother

 \IEEEoverridecommandlockouts 

  \newtheorem{thm}{\bf Theorem}

     \usepackage{stfloats}
  \newtheorem{prop}{\bf Proposition}
  \usepackage{dirtytalk}
\usepackage[]{amsmath}
\usepackage{framed}
\usepackage{amsfonts}
\usepackage{amssymb}
\usepackage{bbm}
\usepackage{breqn}
\usepackage{subcaption}
\usepackage{cite}
\usepackage{enumerate}
\usepackage{cases}
\usepackage{tabu}
\usepackage{url} 
\usepackage{multicol}
\usepackage{color}
\usepackage[bottom]{footmisc}
\usepackage{environ}
\usepackage{tikz} 
\usetikzlibrary{automata,arrows,positioning,calc}
\usepackage[]{algpseudocode}
\algtext*{EndWhile}
\algtext*{EndIf}
\algtext*{EndFor}
\usepackage{algorithm}
\title{\LARGE \bf
Entropy-Regularized Stochastic Games}

\newcommand{\StaticObstacle}[2]{ \fill[black] (#1,#2) rectangle (#1+1,#2+1);}
\newcommand{\initialstate}[2]{ \fill[brown] (#1+0.15,#2+0.15) rectangle (#1+0.9,#2+0.9);}
\newcommand{\initialadversary}[2]{ \fill[red] (#1+0.15,#2+0.15) rectangle (#1+0.9,#2+0.9);}
\newcommand{\goalstate}[2]{ \fill[green] (#1.02,#2.02) rectangle (#1+0.92,#2+0.97);}

\algtext*{EndWhile}
\algtext*{EndIf}
\algtext*{EndFor}
\usepackage{algorithm}
\IEEEoverridecommandlockouts
\author{ Yagiz Savas, Mohamadreza Ahmadi, Takashi Tanaka, Ufuk Topcu  \thanks{ Y. Savas, T. Tanaka and U. Topcu are with the Department of Aerospace Engineering, University of Texas at Austin, TX, USA. E-mail: \{yagiz.savas, ttanaka, utopcu\}@utexas.edu   } \thanks{M. Ahmadi is with the Center for Autonomous Systems and Technologies, California Institute of Technology, CA, USA. E-mail: mrahmadi@caltech.edu } }
\begin{document}
\maketitle
\begin{abstract}  
In two-player zero-sum stochastic games, where two competing players make decisions under uncertainty, a pair of optimal strategies is traditionally described by Nash equilibrium and computed under the assumption that the players have perfect information about the stochastic transition model of the environment. However, implementing such strategies may make the players vulnerable to unforeseen changes in the environment. In this paper, we introduce entropy-regularized stochastic games where each player aims to maximize the causal entropy of its strategy in addition to its expected payoff. The regularization term balances each player's rationality with its belief about the level of misinformation about the transition model. We consider both entropy-regularized $N$-stage and entropy-regularized discounted stochastic games, and establish the existence of a value in both games. Moreover, we prove the sufficiency of Markovian and stationary mixed strategies to attain the value, respectively, in $N$-stage and discounted games. Finally, we present algorithms, which are based on convex optimization problems, to compute the optimal strategies. In a numerical example, we demonstrate the proposed method on a motion planning scenario and illustrate the effect of the regularization term on the expected payoff. 
\end{abstract}
\section{Introduction}
A two-player zero-sum stochastic game (SG) \cite{Shapley} models sequential decision-making of two players with opposing objectives in a stochastic environment. An SG is played in stages. At each stage, the game is in a state, and the players choose one of their available actions simultaneously and receive payoffs. The game then transitions to a new random state according to a probability distribution which represents the stochasticity in the environment.

In an SG, each player aims to synthesize a strategy that maximizes the player's expected payoff at the end of the game. Traditionally, a pair of optimal strategies is described by Nash equilibrium \cite{Nash} according to which both players play their best-response strategies against the opponent's strategy. The value of the game then corresponds to the expected payoff that each player receives at the end of the game, if they both play their respective equilibrium strategies.

The concept of Nash equilibrium is based on the assumptions that the players have perfect information about the environment and act rationally \cite{QRE_book}. However, in certain scenarios, the information that a player has about the environment may not match the reality. For example, in a planning scenario, if the player obtains its information about the environment through surveillance missions performed in the past, it may face with a significantly different environment during the execution of the play. In such scenarios, playing an equilibrium strategy may dramatically decrease the player's actual expected payoff as the strategy is computed under the assumption of perfect information. 

The principle of maximum entropy prescribes a probability distribution that is \say{maximally noncommittal with regard to missing information} \cite{Jaynes}. The principle of maximum causal entropy extends the maximum entropy principle to settings where there is dynamically revealed side information that causally affects the evolution of a stochastic process \cite{Ziebart2, Ziebart3}. A distribution that maximizes the causal entropy of a stochastic process (in the absence of additional constraints) is the one that makes all admissible realizations equally probable regardless of the revealed information \cite{Yagiz}. Therefore, the causal entropy of a player's strategy provides a convenient way to quantify the dependence of its strategy to its level of information about the environment as well as the other player's strategy.

In this paper, we propose a method to synthesize a pair of strategies that balances each player's rationality with its belief about the level of missing information. Specifically, we regularize each player's objective with the causal entropy of its strategy which is causally dependent on the history of play. Therefore, the proposed method allows the players to adjust their strategies according to different levels of misinformation by tuning a parameter that controls the importance of the regularization term. For example, in two extremes, it allows the player to be perfectly rational or to purely randomize its strategy.

We study both entropy-regularized $N$-stage and entropy-regularized discounted games, and show the existence of a value in both games. We first prove the sufficiency of Markovian and stationary strategies for both players, respectively, in $N$-stage and discounted games in order to maximize their entropy-regularized expected payoff. Then, we provide algorithms based on a sequence of convex optimization problems to compute a pair of equilibrium strategies. Finally, we demonstrate the proposed methods on a motion planning scenario, and illustrate that the introduced regularization term yields strategies that perform well in different environments.

\noindent \textbf{Related work.} In stochastic games literature, the idea of balancing the expected payoffs with an additional regularization term appeared recently in \cite{Jordi} and \cite{Reza}. The work \cite{Jordi} proposes to bound the rationality of the players to obtain tunable behavior in video games. They study $\gamma$-discounted games and restrict their attention to stationary strategies to balance the expected payoffs with the Kullback-Leibner distance of the player's strategies from reference strategies. In \cite{Reza}, authors study $N$-stage games and consider only Markovian strategies to balance the expected payoffs with the player's sensing costs which are expressed as directed information from states to actions. Unlike this work, we introduce causal entropy of strategies as the regularization term. Additionally, we allow the player's to follow history-dependent strategies and prove the sufficiency of Markovian strategies to attain the value in $N$-stage games. 

Regularization terms are also used in matrix and extensive form games generally to learn equilibrium strategies \cite{Merti}, \cite{Leslie}, \cite{Ling}, \cite{Levine}. When each player uses the same parameter to regularize its expected payoffs with the entropy of its strategy, an equilibrium strategy profile is called a quantal response equilibrium (QRE) \cite{QRE_book}, and an equilibrium strategy of a player is referred as quantal best response \cite{Leslie} or logit choice strategy \cite{Merti}. From a theoretical perspective, the main difference between our approach and the well-studied QRE concept \cite{McKelvey} is that we establish the existence of equilibrium strategies even if the players use different regularization parameters. Additionally, we provide an efficient algorithm based on a convex optimization problem to compute the equilibrium strategies.

Robust stochastic games \cite{Kardes}, \cite{Aghassi} concern the synthesis of equilibrium strategies when the uncertainty in transition probabilities and payoff functions can be represented by structured sets. Unlike robust SG models, the proposed method in this paper can still be used when it is not possible to form a structured uncertainty set.

In reinforcement learning literature, the use of regularization terms is extensively studied to obtain robust behaviors \cite{Haarnoja}, improve the convergence rates \cite{Fox}, and compute optimal strategies efficiently \cite{Todorov}. As stochastic games model multi-player interactions, our approach leverages the ideas discussed in aforementioned work to environments where an adversary aims to prevent a player to achieve its objective. 

\section{Background}
We first review some concepts from game theory and information theory that will be used in the subsequent sections.

\textbf{Notation:} For a sequence $x$, we write $x^t$ to denote $(x_1,x_2,\ldots, x_t)$. Upper case symbols such as $X$ denote random variables, and lower case symbols such as $x$ denote a specific realization. The cardinality of a set $\mathcal{X}$ is denoted by $\lvert \mathcal{X} \rvert$, and the probability simplex defined over the set $\mathcal{X}$ is denoted by $\Delta(\mathcal{X})$. For $V_1$$,$$V_2$$\in$$\mathbb{R}^n$, we write $V_1$$\preccurlyeq$$V_2$ to denote the coordinate-wise inequalities. We use the index set $\mathbb{Z_{+}}$$=$$\{1,2,\ldots\}$ and the natural logarithm $\log (\cdot)$$=$$\log_e (\cdot)$.
\subsection{Two-Player Stochastic Games}
A two-player stochastic game $\Gamma$ \cite{Shapley} is played in stages. At each stage $t$, the game is in one of its finitely many states $\mathcal{X}$, and each player observes the current state $x_t$. At each state $x_t$, the players choose one of their finitely many actions, and the game transitions to a successor state $x_{t+1}$ according to a probability distribution $\mathcal{P}$$:$$\mathcal{X}$$\times$$\mathcal{U}$$\times$$\mathcal{W}$$\rightarrow$$\Delta(\mathcal{X})$ where $\mathcal{U}$ and $\mathcal{W}$ are finite action spaces for player 1 and player 2, respectively. The pair of actions, $u_t$$\in$$\mathcal{U}$ and $w_t$$\in$$\mathcal{W}$, together with the current state $x_t$$\in$$\mathcal{X}$ determine the payoff $\mathcal{R}(x_t,u_t,w_t)$$\leq$$\overline{\mathcal{R}}$$<$$\infty$ to be made by player 2 to player 1 at stage $t$.

A player's strategy is a specification of a probability distribution over available actions at each stage conditional on the history of the game up to that stage. Formally, let $\mathcal{H}_t$$=$$(\mathcal{X}$$\times$$\mathcal{U}$$\times$$\mathcal{W})^{t-1}$$\times$$\mathcal{X}$ be the set of all possible history of plays up to stage $t$. Then, the strategy of player 1 and player 2 are denoted by $\boldsymbol{\sigma}$$=$$(\sigma_1,\sigma_2,\ldots)$ and $\boldsymbol{\tau}$$=$$(\tau_1, \tau_2,\ldots)$, respectively, where $\sigma_t$ $:$$\mathcal{H}_t$$\rightarrow$$\Delta(\mathcal{U})$ and $\tau_t$ $:$$\mathcal{H}_t$$\rightarrow$$\Delta(\mathcal{W})$ for all $t$. If a player's strategy depends only on the current state for all stages, e.g., $\sigma_t$$:$$\mathcal{X}_t$$\rightarrow$$\Delta(\mathcal{U})$ for all $t$, the strategy is said to be \textit{Markovian}. A \textit{stationary} strategy depends only on the current state and is independent of the stage number, e.g., $\boldsymbol{\sigma}$$=$$(\sigma, \sigma,\ldots)$, where $\sigma$$:$$\mathcal{X}$$\rightarrow$$\Delta(\mathcal{U})$. We denote the set of all strategies, all Markovian strategies, and all stationary strategies for player $i$$\in$$\{1,2\}$ by $\Gamma_i$, $\Gamma_i^M$, and $\Gamma_i^S$, respectively.

Let $\mu_{t+1}(x^{t+1},u^t,w^t)$ be the joint probability distribution over the history $\mathcal{H}_{t+1}$ of play which is uniquely determined by the initial state distribution $\mu_1(h_1)$ through the recursive formula
\begin{align}
&\mu_{t+1}(x^{t+1},u^t,w^t)=\mathcal{P}(x_{t+1} | x_t, u_t, w_t)\sigma_t(u_t | h_t)\nonumber \\
&\qquad\qquad \qquad\qquad \ \times \tau_t(w_t | h_t)\mu_{t}(h_t)
\end{align}
where $h_t$$\in$$\mathcal{H}_t$ is the history of play up to stage $t$.

A stochastic game with the initial distribution $\mu_1(x_1)$ is called an $N$-stage game, if the game ends after $N$ stages. The evaluation function for an $N$-stage game is 
 \begin{align}\label{usual_evaluation_finite}
&J(X^N, U^N, W^N):=\nonumber\\
&\qquad\qquad\sum_{t=1}^N\mathbb{E}^{\overline{\mu}_t}\mathcal{R}(X_t, U_t, W_t)+\mathbb{E}^{{\mu}_{T+1}}\mathcal{R}(X_{N+1})
\end{align}
where $\overline{\mu}_t(\cdot)$$:=$${\mu}_t(\cdot)$$\sigma_t(\cdot)$$\tau_t(\cdot)$. Similarly, if the number of stages in the game is infinite, and the future payoffs are discounted by a factor $0$$<$$\gamma$$<$$1$, the game is called a $\gamma$-discounted game. The evaluation function for a $\gamma$-discounted game is 
\begin{align}\label{classical_infinite}
\sum_{t=1}^{\infty}\gamma^{t-1}\mathbb{E}^{\overline{\mu}_t}\mathcal{R}(X_t, U_t, W_t).
\end{align}

The player 1's objective is to maximize the evaluation function, i.e., its expected payoff, whereas the player 2 aims to minimize it. A stochastic game is said to have the \textit{value} $\mathcal{V}^{\star}$, if for an evaluation function $f(\boldsymbol{\sigma},\boldsymbol{\tau})$, we have
\begin{align*}
\mathcal{V}^{\star}=\max_{\boldsymbol{\sigma}\in \Gamma_u} \min_{\boldsymbol{\tau}\in\Gamma_w}  f(\boldsymbol{\sigma},\boldsymbol{\tau})=\min_{\boldsymbol{\tau}\in\Gamma_w}\max_{\boldsymbol{\sigma}\in \Gamma_u}   f(\boldsymbol{\sigma},\boldsymbol{\tau}).
\end{align*}
A pair of strategies $(\boldsymbol{\sigma}^{\star},\boldsymbol{\tau}^{\star})$ is said to be equilibrium strategies if it attains the value of the game.

It is well-known that both $N$-stage and $\gamma$-discounted games have a value for finite state and action sets \cite{Bewley}. Moreover, Markovian and stationary strategies are sufficient for players to attain the value in $N$-stage and $\gamma$-discounted games, respectively \cite{Shapley}, \cite{Sorin}.

\subsection{Causal Entropy}
For a sequential decision-making problem where decisions depend causally on the past information such as the history of play, the causal entropy of a strategy is a measure to quantify the randomness of the strategy. Let $X^N$, $Y^N$ and $Z^N$ be sequences of random variables with length $N$. The entropy of the sequence $X^N$ causally conditioned on the sequences $Y^N$ and $Z^N$ is defined as \cite{Kramer}
 \begin{align}\label{causal_deff}
H(X^N || Y^N, Z^N ):= \sum_{t=1}^NH(X_t | X^{t-1}, Y^{t}, Z^ t),
\end{align}
where 
\begin{align}
&H(X_t | X^{t-1}, Y^{t}, Z^ t):=\nonumber\\
&-\sum_{\mathcal{X}^N, \mathcal{Y}^{N},\mathcal{Z}^{N} }\text{Pr}(x^t, y^t, z^t)\log \text{Pr}(x_t | x^{t-1}, y^t, z^t).
\end{align}

The concept of causal entropy has recently been used to infer correlated-equilibrium strategies in Markov games \cite{Ziebart1} and to recover cost functions in inverse optimal control problems \cite{Ziebart3}. In this study, we employ causal entropy to compute an equilibrium strategy profile that balances the players' expected payoff with the randomness of their strategies in stochastic games. 

 In the absence of additional constraints, a strategy $\boldsymbol{\sigma}$$\in$$\Gamma_1$ that maximizes the causal entropy $H(U^N || X^N, W^{N-1})$ of the player 1, which is conditioned on the revealed history of play, is the stationary strategy $\boldsymbol{\sigma}$$=$$(\sigma,\sigma,\ldots)$ where $\sigma(x)(u)$$=$$1/ \lvert \mathcal{U} \rvert$. Therefore, a player that maximizes the entropy of its strategy acts purely randomly regardless of the history of play. On the other hand, a player that regularizes its expected payoff with the entropy of its strategy can be thought as a player that balances its rationality with its belief about the correctness of the underlying transition model of the environment.

\section{Problem Statement}
We first consider entropy-regularized $N$-stage games for which we define the evaluation function as
\begin{align}\label{entropy_regularized}
\Phi_N(\boldsymbol{\sigma},\boldsymbol{\tau}) :=& J(X^N, U^N, W^N)+\frac{1}{\beta_1}H(U^N || X^N, W^{N-1} )\nonumber\\
&-\frac{1}{\beta_2}H(W^N || X^N, U^{N-1} ),
\end{align}
where $\beta_1,\beta_2$$>$$0$ are regularization parameters that adjust for players the importance of the randomness in their strategies. 
Note that, when $\beta_1$$=$$\beta_2$$=$$\infty$, both players act perfectly rational, and we recover the evaluation function \eqref{usual_evaluation_finite}. Additionally, since the play is simultaneous, the information of a player's strategy at a given stage is not revealed to the other player. Hence, at each stage, players are allowed to condition their strategies only to observed history of play.

\noindent \textbf{Problem 1:} Provide an algorithm to synthesize, if exists, equilibrium strategies in entropy-regularized $N$-stage games.

We next consider stochastic games that are played in infinite stages, and introduce entropy-regularized $\gamma$-discounted games for which we define the evaluation function as
\begin{align}\label{entropy_regularized_infinite}
\Phi_{\infty}(\boldsymbol{\sigma},\boldsymbol{\tau}):=& \sum_{t=1}^{\infty}\gamma^{t-1}\Big[\mathbb{E}^{\overline{\mu}_t}\mathcal{R}(X_t, U_t, W_t)\nonumber\\
&+\frac{1}{\beta_1}H(U_t | H_t )-\frac{1}{\beta_2}H(W_t | H_t )\Big],
\end{align}
where $H_t$$=$$(X^t, U^{t-1}, W^{t-1})$, i.e., the admissible histories of play at stage $t$. Note that in the evaluation function \eqref{entropy_regularized_infinite}, we discount players' future entropy gains as well as the expected payoff in order to ensure the finiteness of the evaluation function. 

\noindent \textbf{Problem 2:} Provide an algorithm to synthesize, if exists, equilibrium strategies in entropy-regularized $\gamma$-discounted games.

\section{Existence of Values and The Computation of Optimal Strategies}
{\setlength{\parindent}{0cm}
In this section, we analyze entropy regularized $N$-stage and $\gamma$-discounted games, and show that both games have values. Then, we provide algorithms to synthesize equilibrium strategies that attain the corresponding game values.
\subsection{Entropy-Regularized $N$-Stage Games}\label{finite_section}
Searching optimal strategies that solve a stochastic game with the evaluation function $\Phi_N(\boldsymbol{\sigma},\boldsymbol{\tau})$ in the space of all strategies can be intractable for large $N$. We begin with establishing the existence of optimal strategies for both players in the space of Markovian strategies.
\begin{prop}\label{Markovian_finite}
Markovian strategies are sufficient for both players to attain, if exists, the value in entropy-regularized N-stage games, i.e.,
\begin{align*}
\max_{\boldsymbol{\sigma}\in\Gamma_u^M}\min_{\boldsymbol{\tau}\in\Gamma_w^M}\Phi_N(\boldsymbol{\sigma},\boldsymbol{\tau})=\max_{\boldsymbol{\sigma}\in\Gamma_u}\min_{\boldsymbol{\tau}\in\Gamma_w}\Phi_N(\boldsymbol{\sigma},\boldsymbol{\tau}),\\
\min_{\boldsymbol{\tau}\in\Gamma_w^M}\max_{\boldsymbol{\sigma}\in\Gamma_u^M}\Phi_N(\boldsymbol{\sigma},\boldsymbol{\tau})=\min_{\boldsymbol{\tau}\in\Gamma_w}\max_{\boldsymbol{\sigma}\in\Gamma_u}\Phi_N(\boldsymbol{\sigma},\boldsymbol{\tau}).
\end{align*}
\end{prop}}
\noindent\textbf{Proof:} See Appendix A. $\quad \Box$

Next, we show that entropy-regularized $N$-stage games have a value. Let $\rho_t$$:$$\mathcal{X}$$\times$$\mathcal{U}$$\times$$\mathcal{W}$$\rightarrow$$\mathbb{R}$ be a function and $x_t$ be a fixed state. Additionally, let  
\begin{align}\label{one-shot-game}
\mathcal{V}_t^{\sigma_t,\tau_t}(x_t):=\mathbb{E}^{\sigma_t, \tau_t}\Big[&\rho_t(x_t,u_t,w_t)-\frac{1}{\beta_1}\log \sigma_t(u_t | x_t )\nonumber\\
&+\frac{1}{\beta_2}\log \tau_t(w_t | x_t )\Big]
\end{align}
be the evaluation function for a \say{one-shot} game in which the game starts from the state $x_t$ and ends after both players play their one-step strategy. 
{\setlength{\parindent}{0cm}
\begin{prop}\label{normal_game_prop}
A stochastic game with the evaluation function \eqref{one-shot-game} has a value, i.e.,
\begin{align*}
\max_{\sigma_t\in \Delta(\mathcal{U})}\min_{\tau_t\in\Delta(\mathcal{W})}\mathcal{V}_t^{\sigma_t,\tau_t}(x_t)=\min_{\tau_t\in\Delta(\mathcal{W})}\max_{\sigma_t\in \Delta(\mathcal{U})}\mathcal{V}_t^{\sigma_t,\tau_t}(x_t).
\end{align*}
\end{prop}}
\noindent\textbf{Proof:} It is clear that $\mathcal{V}_t(x_t)$ is a continuous function that is concave in $\sigma_t$ and convex in $\tau_t$. Additionally, $\Delta(\mathcal{U})$ and $\Delta(\mathcal{W})$ are compact convex sets. The result follows from von Neumann's minimax theorem \cite{Neumann}.$\quad\Box$

The following proposition states that one can compute the value of the one shot game \eqref{one-shot-game} and synthesize equilibrium strategies by solving a convex optimization problem.

{\setlength{\parindent}{0cm}
\begin{prop}\label{optimal_strategy_finite}
For a given one-shot game with the evaluation function \eqref{one-shot-game}, optimal strategies $(\sigma_t^{\star},\tau_t^{\star})$ satisfy
\begin{align}\label{player_1_st}
&\sigma_t^{\star}(u_t|x_t)\in\nonumber\\
&\arg\max_{\sigma_t\in\Delta(U)}\Big[-\frac{1}{\beta_1}\sum_{u_t\in\mathcal{U}}\sigma_t(u_t | x_t)\log\sigma_t(u_t | x_t)\nonumber\\
&-\frac{1}{\beta_2}\log\sum_{w_t\in \mathcal{W}}\exp\big(-\beta_2\sum_{u_t\in\mathcal{U}}\sigma_t(u_t | x_t)\rho_t(x_t,u_t,w_t)\big)\Big],\\
&\tau_t^{\star}(w_t|x_t)=\nonumber\\\label{player_2_st}
&\frac{\exp\big(-\beta_2\sum_{u_t\in\mathcal{U}}\sigma_t^{\star}(u_t | x_t)\rho_t(x_t,u_t,w_t)\big)}{\sum_{w_t\in \mathcal{W}}\exp\big(-\beta_2\sum_{u_t\in\mathcal{U}}\sigma^{\star}_t(u_t | x_t)\rho_t(x_t,u_t,w_t)\big)}.
\end{align}
Furthermore, the unique value $\mathcal{V}^{\star}_t(x_t)$ of the game is given by
\begin{align}\label{value_game_opt}
&\mathcal{V}_t^{\star}(x_t)=-\frac{1}{\beta_1}\sum_{u_t\in\mathcal{U}}\sigma^{\star}_t(u_t | x_t)\log\sigma^{\star}_t(u_t | x_t)\nonumber\\
&-\frac{1}{\beta_2}\log\sum_{w_t\in \mathcal{W}}\exp\big(-\beta_2\sum_{u_t\in\mathcal{U}}\sigma^{\star}_t(u_t | x_t)\rho_t(x_t,u_t,w_t)\big).
\end{align}
\end{prop}}
\noindent\textbf{Proof:} See Appendix A.$\quad \Box$

It is worth noting that the objective function of the optimization problem given in \eqref{player_1_st} is strictly concave, and therefore, optimal strategies $\sigma_t^{\star}(u_t|x_t)$ and $\tau_t^{\star}(w_t|x_t)$ are unique. Additionaly, an optimal strategy with the form \eqref{player_2_st} is known in the economics literature as quantal best response \cite{McKelvey}, and for $\beta_1$$=$$\beta_2$$<$$\infty$, the optimal strategies form the well-studied quantal response equilibrium strategies \cite{QRE_book}.

We remark that the optimization problem in \eqref{player_1_st} has a closed-form solution which is a function of the optimal strategy $\tau_t^{\star}(w_t|x_t)$. However, since the closed-form expressions for equilibrium strategies constitute a system of coupled nonlinear equations, the convex optimization formulation provides a more convenient way to compute equilibrium strategies.

Utilizing the results of above propositions, we now reformulate an entropy-regularized $N$-stage game as a series of \say{one-shot} games through the use of Bellman recursions. Let 
\begin{align}\label{recursion_function}
&\rho_t(x_t,u_t,w_t)=\mathcal{R}(x_t,u_t,w_t)\nonumber\\
&\qquad+\sum_{x_{t+1}\in \mathcal{X}}\mathcal{P}(x_{t+1}| x_t, u_t, w_t)\mathcal{V}_{t+1}^{\sigma_t,\tau_t}(x_{t+1}),
\end{align}
for $t$$=$$1,\ldots, N$ where $\mathcal{V}_{N+1}^{\sigma_t,\tau_t}(x_{N+1})$$=$$\mathcal{R}(X_{N+1})$. Then, it can be easily verified that 
\begin{align}\label{objective-equivalence}
\Phi_N(\boldsymbol{\sigma},\boldsymbol{\tau})=\sum_{x_1\in\mathcal{X}}\mu_1(x_1)\mathcal{V}_1(x_1)
\end{align}
for a given initial distribution $\mu_1(x_1)$. Consequently, we obtain the following result.
{\setlength{\parindent}{0cm}
\begin{thm}\label{finite_game_has_value}
Entropy-regularized $N$-stage games have a value. 
\end{thm}}
\noindent\textbf{Proof:} Due to Proposition \ref{Markovian_finite}, we can focus on Markovian strategies to find an equilibrium point in $N$-stage games. We start from the stage $k$$=$$N$ and compute the value of the one shot game \eqref{one-shot-game}, which exists due to Proposition \ref{normal_game_prop}. Using \eqref{one-shot-game} with \eqref{recursion_function} for $k$$=$$N$$-$$1,N$$-$$2,\ldots,1$, we compute the value of $N$$-$$k$$+$$1$ stage games. As a result, the claim follows due to the equivalence given in \eqref{objective-equivalence}.$\quad\Box$ 

Algorithm \ref{euclid} summarizes the computation of the pair $(\boldsymbol{\sigma}^{\star},\boldsymbol{\tau}^{\star})$ of optimal strategies for entropy-regularized $N$-stage games.

\begin{algorithm}
\caption{Strategy computation for $N$-stage games}\label{euclid}
\begin{algorithmic}[1]

\State \textbf{Initialize:} $\mathcal{V}_{N+1}(x_{N+1})$=$\mathcal{R}(x_{N+1})$ for all $x_{N+1}$$\in$$\mathcal{X}$.
      \For{$t = N, N$$-$$1, ..., 1$}
      \State Compute $\rho_t(x_t,u_t,w_t)$ for all $x_t$$\in$$\mathcal{X}$, $u_t$$\in$$\mathcal{U}$, and $w_t$$\in$$\mathcal{W}$ as in \eqref{recursion_function}.
        \State For all $x_t$$\in$$\mathcal{X}$, \begin{align*} &\mathcal{V}^{\star}_{t}(x_{t})=\max_{\sigma_t\in\Delta(U)}-\frac{1}{\beta_1}\sum_{u_t\in\mathcal{U}}\sigma_t(u_t | x_t)\log\sigma_t(u_t | x_t)\\
&-\frac{1}{\beta_2}\log\sum_{w_t\in \mathcal{W}}\exp(-\beta_2\sum_{u_t\in\mathcal{U}}\sigma_t(u_t | x_t)\rho_t(x_t,u_t,w_t))\end{align*}
\State For all $x_t$$\in$$\mathcal{X}$, compute $\sigma_t^{\star}$ and $\tau_t^{\star}$ as in \eqref{player_1_st} and \eqref{player_2_st}, respectively.
      \EndFor

\State \textbf{return}  $\boldsymbol{\sigma}^{\star}$$=$$(\sigma_1^{\star},\ldots, \sigma_T^{\star})$ and $\boldsymbol{\tau}^{\star}$$=$$(\tau_1^{\star},\ldots, \tau_T^{\star})$.
\end{algorithmic}
\end{algorithm}
\noindent\textbf{Remark:} In certain scenarios, one of the players may prefer to play perfectly rationally against a boundedly rational opponent, e.g., $\beta_2$$=$$\infty$. In that case, it is still possible to compute equilibrium strategies by solving a convex optimization problem at each stage. The value of the one-shot game \eqref{one-shot-game} still exists due to the arguments provided in the proof of Proposition \ref{normal_game_prop}. However, the form of optimal strategies slightly changes to
\begin{align}\label{player_1_infinity}
&\tau_t^{\star}(w_t|x_t)=\arg\min_{\tau_t\in\Delta(W)}\frac{1}{\beta_1}\log\sum_{u_t\in \mathcal{U}}\exp\Big(\nonumber\\
&\qquad\qquad\qquad\beta_1\sum_{w_t\in\mathcal{W}}\tau_t(w_t | x_t)\rho_t(x_t,u_t,w_t)\Big),\\
&\sigma_t^{\star}(u_t|x_t)=\nonumber\\\label{player_2_infinity}
&\frac{\exp\big(\beta_1\sum_{w_t\in\mathcal{W}}\tau_t^{\star}(u_t | x_t)\rho_t(x_t,u_t,w_t)\big)}{\sum_{u_t\in \mathcal{U}}\exp\big(\beta_1\sum_{w_t\in\mathcal{W}}\tau^{\star}_t(w_t | x_t)\rho_t(x_t,u_t,w_t)\big)}.
\end{align}
It is important to note that, if $\beta_i$$=$$\infty$ for some $i$$\in$$\{1,2\}$, equilibrium strategies may be not unique since the function $\log\sum\exp(\cdot)$ is not strictly convex over its domain \cite{Rockafellar}.

\subsection{Entropy-Regularized $\gamma$-Discounted Games}
In this section, we focus on Markovian strategies, whose optimality for $N$-stage games is shown in Proposition \ref{Markovian_finite}. Let $\mathcal{V}$$\in$$\mathbb{R}^{\lvert \mathcal{X}\rvert}$ be a real-valued function, and for a given $x$$\in$$\mathcal{X}$, $\mathcal{L}(\mathcal{V})(x,\cdot,\cdot)$ $:$ $\mathcal{U}$$\times$$\mathcal{W}$$\rightarrow$$\mathbb{R}$ be a function such that
\begin{align}
&\mathcal{L}(\mathcal{V})(x,\sigma,\tau):=\mathbb{E}^{\sigma, \tau}\Big[\mathcal{R}(x,u,w)-\frac{1}{\beta_1}\log \sigma(u | x )\nonumber\\
&+\frac{1}{\beta_2}\log \tau(w | x )+\gamma\sum_{{x'}\in\mathcal{X}}\mathcal{P}(x' | x,u, w)\mathcal{V}(x')\Big].
\end{align}

As discussed in Proposition \ref{normal_game_prop}, a one-shot game with the evaluation function $\mathcal{L}(\mathcal{V})(x,\sigma,\tau)$ has a value. Therefore, we can introduce the Shapley operator $\Psi$$:$$\mathcal{V}$$\rightarrow$$\Psi(\mathcal{V})$ from $\mathbb{R}^{\lvert \mathcal{X}\rvert}$ to itself specified, for all $x$$\in$$\mathcal{X}$, as 
\begin{align*}
\Psi(\mathcal{V})[x]:=\max_{\sigma\in\Delta(U)}\min_{\tau\in\Delta(W)}\mathcal{L}(\mathcal{V})(x,\sigma,\tau).
\end{align*}

It is clear that the operator $\Psi$ satisfies two key properties: \textit{monotonicity}, i.e., $\mathcal{V}$$\preccurlyeq$$\overline{\mathcal{V}}$ implies $\Psi(\mathcal{V})$$\preccurlyeq$$\Psi(\overline{\mathcal{V}})$, and \textit{reduction of constants}, i.e., for any $k$$\geq$$0$, $\Psi(\mathcal{V}$$+$$k\boldsymbol{1})[x]$$=$$\Psi(\mathcal{V})[x]$$+$$\gamma k$ for all $x$$\in$$\mathcal{X}$. Consequently, it is straightforward to show that the operator $\Psi$ is a contraction mapping \cite{Bertsekas}. Specifically, we have
\begin{align*}
\lVert \Psi(\mathcal{V})-\Psi(\overline{\mathcal{V}})\rVert_{\infty}\leq \gamma \lVert \mathcal{V}- \overline{\mathcal{V}}\rVert_{\infty},
\end{align*}
where $\lVert\mathcal{V}\rVert_{\infty}$$=$$\max_{x\in\mathcal{X}}\mathcal{V}(x)$. We omit the details since similar results can be easily found in the literature \cite{Neyman}. Then, by Banach's fixed-point theorem \cite{Puterman}, we conclude that the operator $\Psi$ has a unique fixed point which satisfies $\mathcal{V}^{\star}$$=$$\Psi\mathcal{V}^{\star}$.

Next, we need to show that the fixed point $\mathcal{V}^{\star}$ is indeed the value of the entropy-regularized $\gamma$-discounted game. Let $\boldsymbol{\sigma}^{\star}$$=$$(\sigma^{\star},\sigma^{\star},\ldots)$ be a stationary strategy such that $\sigma^{\star}$ is a one-step strategy for player 1 satisfying the fixed point equation, and $\boldsymbol{\tau}$ be an arbitrary Markovian strategy for player 2. Denoting by $h_t$ the history of play of length $t$, one has, by definition of $\Psi$ and $\boldsymbol{\sigma}^{\star}$, 
\begin{align*}
&\mathbb{E}^{\overline{\mu}_t}\Big[\mathcal{R}(x_t,u_t,w_t)-\frac{1}{\beta_1}\log \sigma^{\star}(u_t | x_t )+\frac{1}{\beta_2}\log \tau(w_t | x_t )\nonumber\\
&+\gamma\sum_{{x'}\in\mathcal{X}}\mathcal{P}(x' | x_t,u_t, w_t)\mathcal{V}^{\star}(x') \Big| h_t\Big]\geq \mathbb{E}^{\overline{\mu}_t}\Big[\mathcal{V}^{\star}(x_t)\Big| h_t\Big].
\end{align*}
This expression can further be written as
\begin{align*}
&\mathbb{E}^{\overline{\mu}_t,\overline{\mu}_{t+1}}\Big[\mathcal{R}(x_t,u_t,w_t)-\frac{1}{\beta_1}\log \sigma^{\star}(u_t | x_t )\nonumber\\
&+\frac{1}{\beta_2}\log \tau(w_t | x_t )+\gamma\mathcal{V}^{\star}(x_{t+1})\Big| h_t\Big]\geq \mathbb{E}^{\overline{\mu}_t}\Big[\mathcal{V}^{\star}(x_t)\Big| h_t\Big].
\end{align*}
Multiplying by $\gamma^{t-1}$, taking expectations and summing over $1$$\leq$$t$$<$$k$, one obtains
\begin{align*}
&\sum_{t=1}^{k-1}\gamma^{t-1}\mathbb{E}^{\overline{\mu}_t}\Big[\mathcal{R}(x_t,u_t,w_t)-\frac{1}{\beta_1}\log \sigma^{\star}(u_t | x_t )+\nonumber\\
&\frac{1}{\beta_2}\log \tau(w_t | x_t )\Big| x_1\Big]\geq \mathcal{V}^{\star}(x_1)-\gamma^{k}\mathbb{E}^{\overline{\mu}_{k+1}}\Big[\mathcal{V}^{\star}(x_{k+1})\Big| x_1\Big].
\end{align*}
Taking the limit as $k$$\rightarrow$$\infty$ and using Proposition \ref{Markovian_finite}, we obtain
\begin{align}\label{infinite_first}
\Phi_{\infty}(\sigma^{\star},\tau)\geq \mathcal{V}^{\star}(x_1).
\end{align}
Similarly, when player 2 plays the optimal stationary strategy $\boldsymbol{\tau}^{\star}$$=$$(\tau^{\star},\tau^{\star},\ldots)$ against an arbitrary Markovian strategy $\boldsymbol{\sigma}$ of player 1, we have 
\begin{align}\label{infinite_second}
\Phi_{\infty}(\sigma,\tau^{\star})\leq \mathcal{V}^{\star}(x_1).
\end{align}
Then, the combination of \eqref{infinite_first} and \eqref{infinite_second} implies the following result.
{\setlength{\parindent}{0cm}
\begin{thm}
Entropy-regularized $\gamma$-discounted games have a value which satisfies $\Psi(\mathcal{V})$$=$$\mathcal{V}$. Furthermore, stationary strategies are sufficient for both players to attain the game value, i.e.,
\begin{align*}
\max_{\boldsymbol{\sigma}\in\Gamma}\min_{\boldsymbol{\tau}\in\Gamma}\Phi_{\infty}(\sigma,\tau)=\max_{\boldsymbol{\sigma}\in\Gamma^S}\min_{\boldsymbol{\tau}\in\Gamma^S}\Phi_{\infty}(\sigma,\tau).
\end{align*}
\end{thm}}
Computation of optimal strategies is just an extension of Algorithm \ref{euclid}. Note that for $\gamma$-discounted games, we use the same one-shot game introduced in \eqref{one-shot-game}.  Therefore, optimal decision rules at each stage has the form \eqref{player_1_st} and \eqref{player_2_st}. Consequently, to compute the optimal strategies, we initialize Algorithm \ref{euclid} with an arbitrary value vector $\mathcal{V}$$\in$$\mathbb{R}^{\lvert \mathcal{X}\rvert}$ and iterate until convergence, which is guaranteed by the existence of a unique fixed point.

\section{A Numerical Example}
In this section, we demonstrate the proposed strategy synthesis method on a motion planning scenario that we model as an entropy-regularized $\gamma$-discounted game. To solve the convex optimization problems required for the computation of equilibrium strategies, we use ECOS solver \cite{Ecos} through the interface of CVXPY \cite{cvxpy}. All computations are performed by setting $\gamma$$=$$0.8$. 

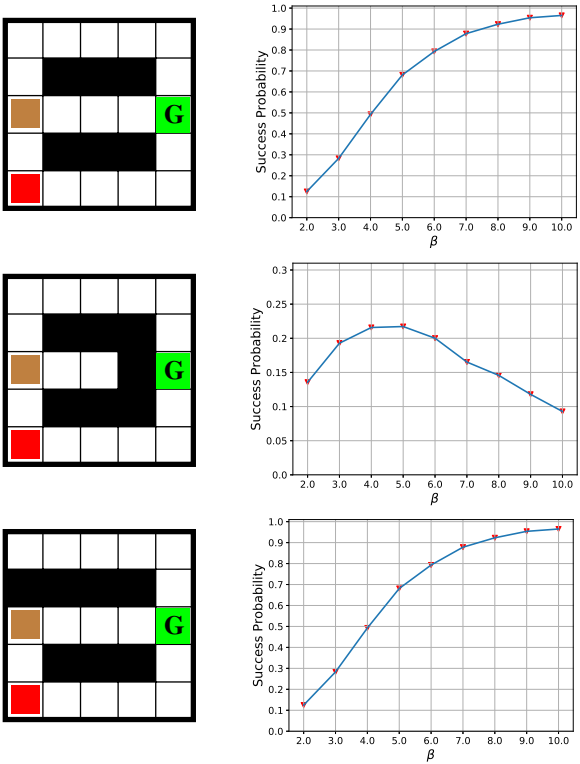
\begin{figure}[t!]\centering\hspace{-0.5\linewidth}
\begin{subfigure}[!t]{0.1\linewidth}\centering
\raisebox{0.07\height}{\scalebox{0.5}{
\begin{tikzpicture}
\draw[black,line width=1pt] (0,0) grid[step=1] (5,5);
\draw[black,line width=4pt] (0,0) rectangle (5,5);
			   \StaticObstacle{1}{1}  \StaticObstacle{2}{1}  \StaticObstacle{3}{1}
			   \StaticObstacle{1}{3}  \StaticObstacle{2}{3}  \StaticObstacle{3}{3}
			    \initialstate{0}{2}  \goalstate{4}{2} \initialadversary{0}{0} 
\node at (4+0.5,2+0.5) { \textbf{\huge G}};
\end{tikzpicture}
}}
\end{subfigure}
\hspace{0.25\linewidth}
\scalebox{0.3}{
\begin{subfigure}[!t]{0.4\linewidth}
\includegraphics[]{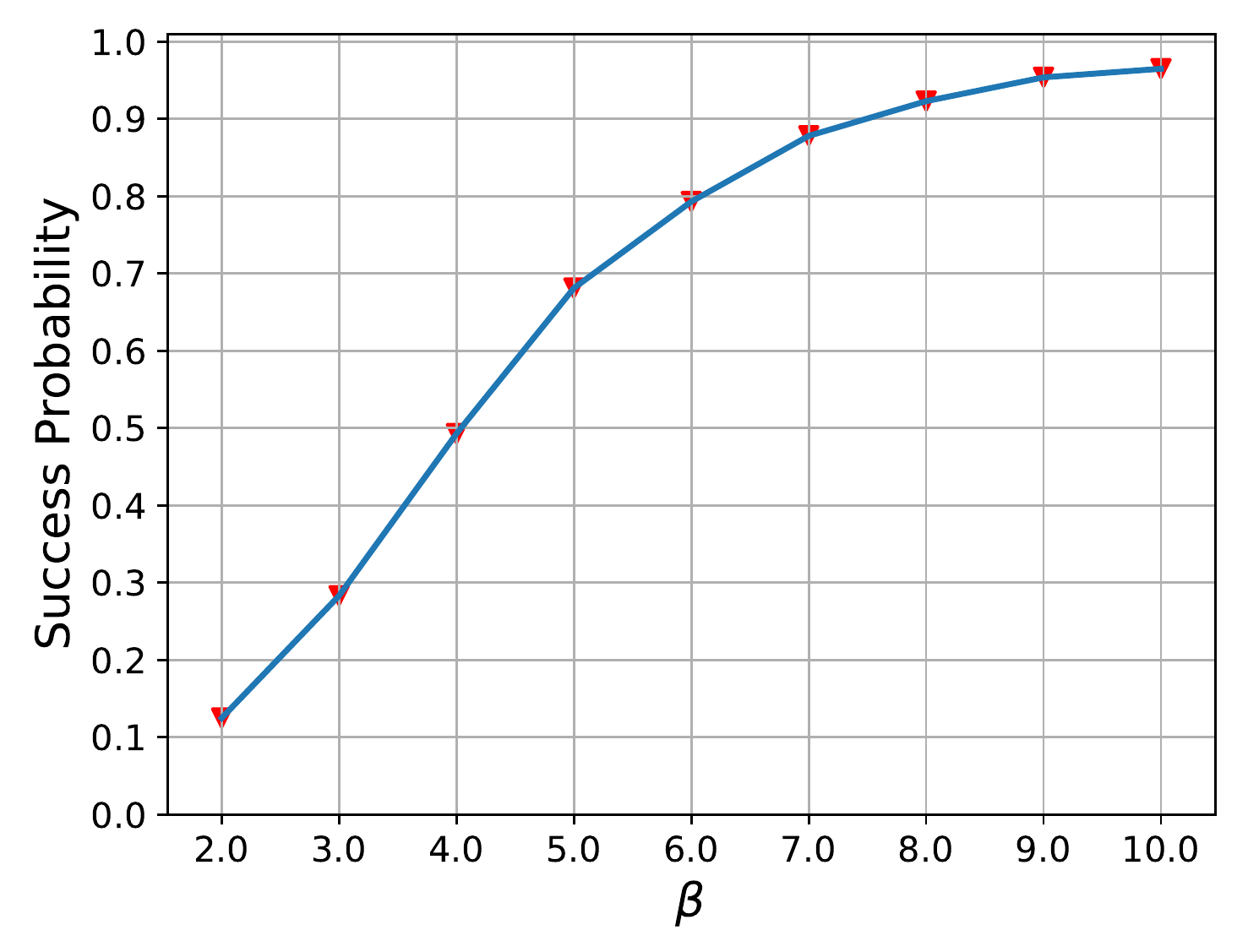}
\end{subfigure}}
\\
\centering\hspace{-0.5\linewidth}
\begin{subfigure}[!t]{0.1\linewidth}\centering
\raisebox{0.07\height}{\scalebox{0.5}{
\begin{tikzpicture}
\draw[black,line width=1pt] (0,0) grid[step=1] (5,5);
\draw[black,line width=4pt] (0,0) rectangle (5,5);
			   \StaticObstacle{1}{1}  \StaticObstacle{2}{1}  \StaticObstacle{3}{1} \StaticObstacle{3}{2}
			   \StaticObstacle{1}{3}  \StaticObstacle{2}{3}  \StaticObstacle{3}{3}
			    \initialstate{0}{2}  \goalstate{4}{2} \initialadversary{0}{0} 
\node at (4+0.5,2+0.5) { \textbf{\huge G}};
\end{tikzpicture}
}}
\end{subfigure}
\hspace{0.25\linewidth}
\scalebox{0.3}{
\begin{subfigure}[!t]{0.4\linewidth}
\includegraphics[]{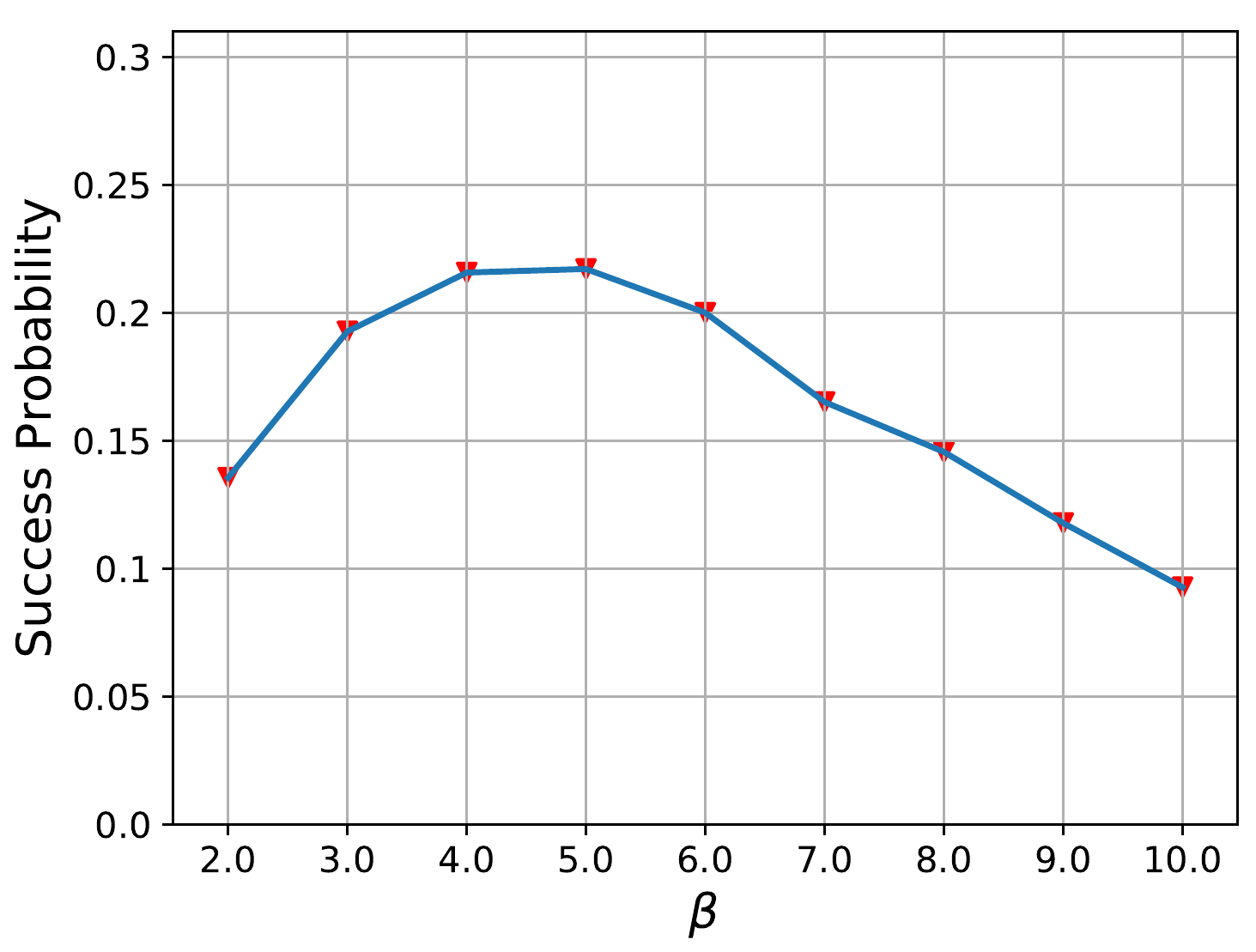}
\end{subfigure}}
\\
\centering\hspace{-0.5\linewidth}
\begin{subfigure}[!t]{0.1\linewidth}\centering
\raisebox{0.07\height}{\scalebox{0.5}{
\begin{tikzpicture}
\draw[black,line width=1pt] (0,0) grid[step=1] (5,5);
\draw[black,line width=4pt] (0,0) rectangle (5,5);
			   \StaticObstacle{1}{1}  \StaticObstacle{2}{1}  \StaticObstacle{3}{1} \StaticObstacle{0}{3}
			   \StaticObstacle{1}{3}  \StaticObstacle{2}{3}  \StaticObstacle{3}{3}
			    \initialstate{0}{2}  \goalstate{4}{2} \initialadversary{0}{0} 
\node at (4+0.5,2+0.5) { \textbf{\huge G}};
\end{tikzpicture}
}}
\end{subfigure}
\hspace{0.25\linewidth}
\scalebox{0.3}{
\begin{subfigure}[!t]{0.4\linewidth}
\includegraphics[]{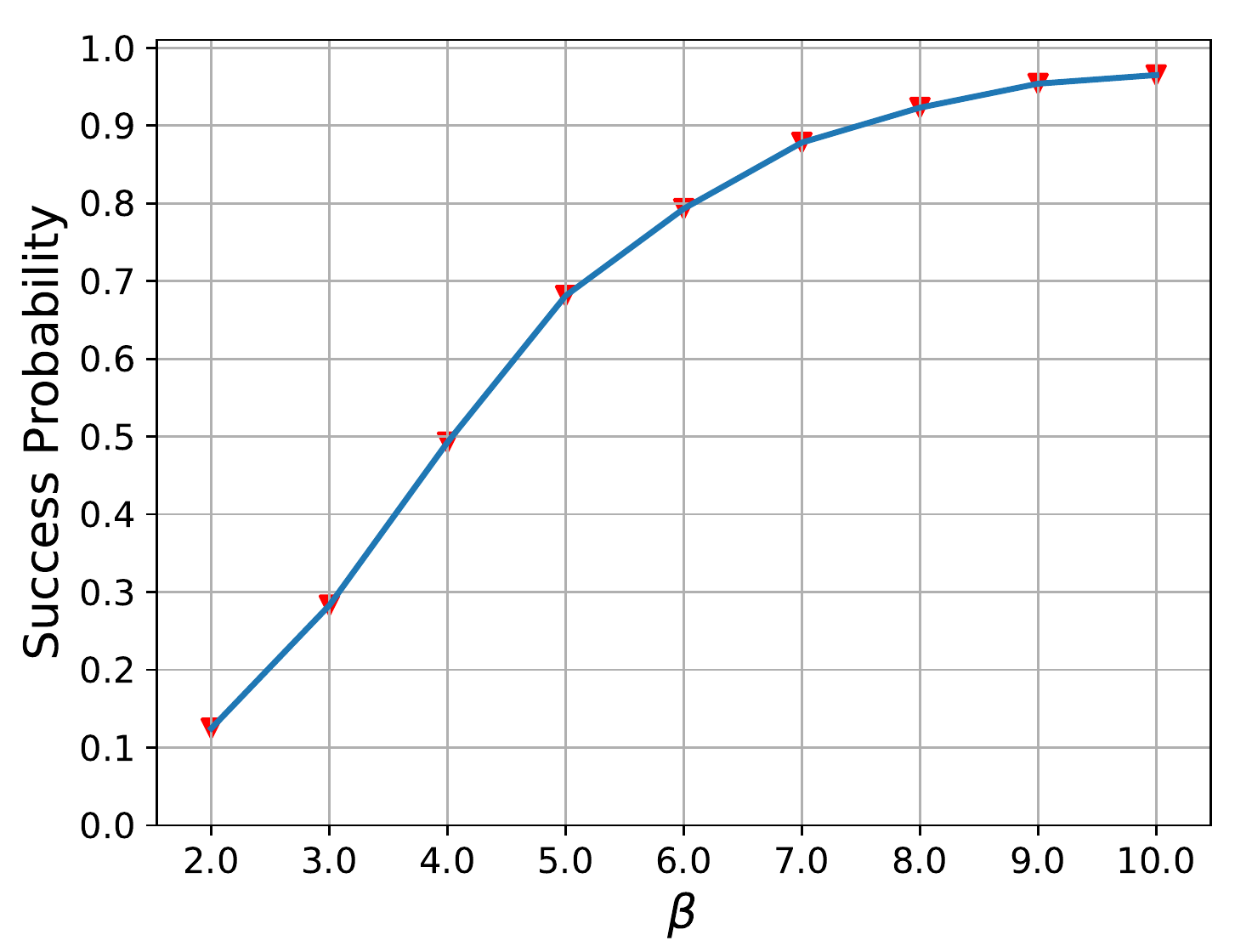}
\end{subfigure}}
\caption{ (Top left) The nominal environment players use for computing their strategies. (Top right) The probability that the player 1 wins the game when it plays the strategy computed by using $\beta_1=\beta_2=\beta$ against the perfectly rational player 2. (Middle and bottom left) The actual environments where the game is played.  (Middle and bottom right) The probability of winning for the player 1 when it employs strategies computed by using different $\beta$ values against the perfectly rational player 2.  }
\label{grid_graph}
\end{figure}

 As the environment model, we consider a $5$$\times$$5$ grid world which is given in Figure \ref{grid_graph} (top left). The brown grid denotes the initial position of the player 1 which aims to reach the goal (green) state. The red grid is the initial position of the player 2 whose aim is to catch the player 1 before reaching the goal state. Finally, black grids represent walls.

Let $x$$=$$(s_1,s_2)$ be the current state of the game such that $x[1]$$=$$s_1$ and $x[2]$$=$$s_2$ are the positions of the player 1 and the player 2, respectively. At each state, the action space for both players is given as $\mathcal{U}$$=$$\mathcal{W}$$=$$\{right,left,up,down,stay\}$. For simplicity, we assume deterministic transitions, i.e., $\mathcal{P}(x,u,w)$$\in$$\{0,1\}$ for all $x$$\in$$\mathcal{X}$, $u$$\in$$\mathcal{U}$ and $w$$\in$$\mathcal{W}$. If a player takes an action for which the successor state is a wall, the player stays in the same state with probability 1.

For a given $(x,u,w)$$\in$$\mathcal{X}$$\times$$\mathcal{U}$$\times$$\mathcal{W}$, we encode the payoff function $\mathcal{R}(x,u,w)$ as the sum of two functions such that $\mathcal{R}(x,u,w)$$=$$\mathcal{R}_1(x,u,w)$$+$$\mathcal{R}_2(x,u,w)$ where 
\begin{align*}
&\mathcal{R}_1(x,u,w)=\sum_{x'[1]=\text{G}}\mathcal{P}(x' | x,u,w), \\
&\mathcal{R}_2(x,u,w)=\begin{cases} -\mathcal{P}(x' | x,u,w) & \text{if} \ \ x'[1]=x'[2]\neq \text{G}\\ 0 & \text{otherwise}.\end{cases}
\end{align*}

Note that the payoff function defines a zero-sum game which is won by the player 1, if it reaches the goal state before getting caught, and by the player 2, if it catches the player 1 before reaching the goal state. 

We first compute Nash equilibrium strategies in the absence of causal entropy terms, i.e., $\beta_1$$=$$\beta_2$$=$$\infty$, by employing standard linear programming formulation \cite{Miltersen} for zero-sum games. Starting from the initial state, an equilibrium strategy for the player 1 is to move towards the goal state by taking the action $right$ deterministically, and for the player 2 is to chase the player 1 by taking the action $up$ in the first two stages, and then, take the action $right$ until reaching the goal state. Therefore, a perfectly rational player 1 wins the game with probability 1 no matter what strategy is followed by the player 2. 

To illustrate the drawback of playing with perfect rationality, we assume that there is another wall in the environment about which the players have no information while they compute the equilibrium strategies, i.e., the players use the nominal environment (top left) to compute the equilibrium strategies. First, we consider the case that the wall is between the goal state and the player 1, as shown in Figure \ref{grid_graph} (middle left). In this case, if the player 1 follows the Nash equilibrium strategy, the probability that it reaches the goal state becomes zero. Therefore, following the Nash equilibrium strategy makes the player 1 significantly vulnerable to such changes in the environment. 

To investigate the tradeoff between rationality and randomness, we compute 9 different strategies for player 1 by using $\beta_1$$=$$\beta_2$$=$$2,3,\ldots,10$, and let it play against the perfectly rational player 2 which follows its Nash equilibrium strategy computed on the nominal environment (top left). The winning probabilities of player 1 under different strategies are shown in Figure \ref{grid_graph} (right) for the corresponding environments given in Figure \ref{grid_graph} (left). This specific scenario demonstrates that, by choosing $\beta$$=$$6$, the player 1 can obtain a robust behavior against unforeseen changes in the environment, i.e., the winning probability is around 20\%, without sacrificing too much from its optimal performance, i.e., around 15\%, if the structure of the environment remains the same. It is worth noting that the asymptotical performance of the player 1 as $\beta$$\rightarrow$$\infty$ approaches to its performance under Nash equilibrium strategy as discussed in \cite{McKelvey}. Additionally, the importance of the randomness for the player 1 increases as $\beta$$\rightarrow$$0$, and using smaller $\beta$ values negatively affects the performance after a critical point, i.e., $\beta$$=$$4$.

Finally, one can argue that the tradeoff occurs in this specific scenario only if the unpredicted wall is between the player 1 and the goal state. To justify the choice of $\beta$ value, we also consider the scenario in which the unexpected wall occupies another state which is shown in Figure \ref{grid_graph} (bottom left). In this case, as shown in  Figure \ref{grid_graph} (bottom right), the use of $\beta$$=$$6$ result in a strategy that guarantees around 80\% winning probability. Therefore, the entropy-regularized strategy of the player 1 still provides an advantage against unpredicted changes without sacrificing too much from the optimal performance.

\section{Conclusions and Future Work}
We consider the problem of two-player zero-sum stochastic games with entropy regularization, wherein the players aim to maximize the causal entropy of their strategies in addition to the conventional expected payoff. We show that equilibrium strategies exist for both entropy-regularized $N$-stage and entropy-regularized $\gamma$-discounted games, and can be computed by solving a convex optimization problem. In numerical examples, we applied the proposed approach to a motion planning scenario and observed that by tuning the regularization parameter, a player can synthesize robust strategies that perform well in different environments against a perfectly rational opponent.

Extending this work to multi-agent reinforcement learning settings, as discussed in \cite{Jordi}, is an interesting future direction. Future work can also investigate the effect of entropy regularization term on the convergence rate of learning algorithms as discussed in \cite{Fox}.

\bibliography{AAAI_2.bib}
\bibliographystyle{IEEEtran}
\section*{Appendix A}
\textbf{Proof of Proposition \ref{Markovian_finite}:} We show the sufficiency of Markovian strategies only for the maximin problem. The proof for the minimax formulation follows the same lines with the arguments provided below. 

The proof is based on backward induction on the stage number $1$$\leq$$k$$\leq$$N$. Let 
\begin{align*}
&\mathcal{V}_k:=\max_{\boldsymbol{\sigma}\in \Gamma_u}\min_{\boldsymbol{\tau}\in\Gamma_w}\sum_{l=k}^N\mathbb{E}^{\overline{\mu}_l}\Big[\mathcal{R}(X_l, U_l, W_l)\nonumber\\
&-\frac{\log \sigma_l(U_l | H_l )}{\beta_1}+\frac{\log \tau_l(W_l | H_l )}{\beta_2}\Big]+\mathbb{E}^{{\mu}_{T+1}}\mathcal{R}(X_{N+1})
\end{align*}
be the value of the $N$$-$$k$ stage problem. Then, we can write the value of $N$$-$$k$ stage problem recursively as
\begin{align}\label{recursion}
&\mathcal{V}_k=\max_{\sigma_k}\min_{\tau_k}\mathbb{E}^{\overline{\mu}_k}\Big[\mathcal{R}(X_k, U_k, W_k)-\frac{1}{\beta_1}\log \sigma_k(U_k | H_k )\nonumber\\
&\qquad \qquad\ +\frac{1}{\beta_2}\log \tau_k(W_k | H_k )+\mathbb{E}^{\mathcal{P}}[\mathcal{V}_{k+1}]\Big].
\end{align}
\underline{Base step: $k$$=$$N$.} Let $\sigma_N$ and $\tau_N^{{\star}}$ be an arbitrary strategy for player 1 and the optimal strategy for player 2 at stage $N$, respectively. Let
\begin{align*}
&\lambda_N(h_N, u_N, w_N):=\mu_N(x^N,u^{N-1},w^{N-1})\\
&\qquad\qquad\qquad\qquad\qquad\qquad\times\sigma_N(u_N| h_N)\tau_N^{{\star}}(w_N| h_N)
\end{align*} 
be the joint distribution induced by $\sigma_N(u_N| h_N)$ and $\tau_N^{\star}(w_N| h_N)$. Additionally, let $\lambda_N(x_N,w_N)$ and $\lambda_N(x_N)$ be the marginal distributions of $\lambda_N(h_N, u_N, w_N)$. We construct a new strategy for player 2 as $\overline{\tau}_N(w_N| x_N)$$:=$$\frac{\lambda_N(x_N,w_N)}{\lambda_N(x_N)}$. Let 
\begin{align*}
&\overline{\lambda}_N(h_N, u_N, w_N):=\mu_N(x^N, u^{N-1},w^{N-1})\\
&\qquad\qquad\qquad\qquad\qquad\qquad\times\sigma_N(u_N| h_N)\overline{\tau}_N(w_N| x_N)
\end{align*}
be the joint distribution induced by $\overline{\tau}_N(w_N| x_N)$. Then, by construction, we have $\overline{\lambda}_N(x_N, u_N, w_N)$$=$${\lambda}_N(x_N, u_N, w_N)$, which can be easily verified by calculating the corresponding marginal distributions. (For a similar strategy construction, see Theorem 5.5.1 in \cite{Puterman}.)

The inner optimization problem in \eqref{recursion} for $k$$=$$N$ reads
\begin{align}
\mathcal{V}_N=\min_{\tau_N}J_N^c(\lambda_N)-J_N^H(\lambda_N)
\end{align}
where 
\begin{align*}
&J_N^c(\lambda_N)=\mathbb{E}^{\lambda_N, \mathcal{P}}[\mathcal{R}(X_N,U_N, W_N)+\mathcal{R}(X_{N+1})\\
&\qquad\qquad\qquad\qquad \qquad-\frac{1}{\beta_1}\log \sigma_N(U_N | H_N)]\\
&J_N^H(\lambda_N)=\frac{1}{\beta_2}H_{\lambda_N}(W_N| X^N, U^{N-1}, W^{N-1}).
\end{align*}
Since the strategy $\sigma_N$ is arbitrarily chosen, it is sufficient to show that $J_N^c(\lambda_N)$$=$$J_N^c(\overline{\lambda}_N)$ and $J_N^H(\lambda_N)$$\leq$$J_N^H(\overline{\lambda}_N)$ in order to establish the sufficiency of Markovian strategies for player 2. The first equality holds by construction. (Note that the $\log(\cdot)$ term is indifferent to changes in the strategy of player 2.) The second inequality can be derived as 
\begin{align}\label{ent_1}
H_{\lambda_N}(W_N| X^N, U^{N-1}, W^{N-1})\leq H_{\lambda_N}(W_N| X_N)&\\\label{ent_2}
= H_{\overline{\lambda}_N}(W_N| X_N)&\\\label{ent_3}
= H_{\overline{\lambda}_N}(W_N|  X^N, U^{N-1}, W^{N-1})&
\end{align}
where \eqref{ent_1} holds since conditioning reduces entropy \cite{Cover}, \eqref{ent_2} is because $\lambda_N$$=$$\overline{\lambda}_N$ by construction, and \eqref{ent_3} is due to the fact that $\overline{\tau}_N(w_N | x_N )$ is a Markovian strategy. 
 Consequently, for any strategy chosen by player 1 in stage $N$, player 2 has a best response strategy in the space of Markovian strategies. 

Next, we can assume that player 2 uses a Markovian strategy and show through a similar strategy construction explained above that player 1 has an optimal strategy in the space of Markovian strategies. As a result, the value $\mathcal{V}_N$ depends on the joint distribution $\mu_N(x^N, u^{N-1}, w^{N-1})$ only through its marginal $\mu_N(x_N)$ and becomes a function of the marginal $\mu_N(x_N)$ only.

\underline{Inductive step: $k$$=$$t$} Assume that Markovian strategies suffice for both players for $k$$=$$t$$+$$1,t$$+$$2,\ldots,N$. Then, by induction hypothesis, $\mathcal{V}_{t+1}$ is a function of $\mu_{t+1}(x_{t+1})$ only. Therefore, using the similar construction to the case $k$$=$$N$, we can construct Markovian strategies $\overline{\sigma}_t(u_t | x_t)$ and $\overline{\tau}_t(w_t | x_t)$ such that the objective function in the right hand side of \eqref{recursion} attained by $\overline{\sigma}_t$ and $\overline{\tau}_t$ is equal to the value of $N-t$ stage problem. As a result, we conclude that Markovian strategies are sufficient for both players to solve the maximin problem \eqref{entropy_regularized}. $\quad \Box$

\textbf{Proof of Proposition \ref{optimal_strategy_finite}:}
Since the one-shot game has a value, without loss of generality, we focus on the problem $\max_{\sigma_t\in \Delta(\mathcal{U})}\min_{\tau_t\in\Delta(\mathcal{W})}\mathcal{V}_t(x_t)$. For notational convenience, we rewrite the problem as
\begin{align*}
&\max_{Q^{ij}}\min_{Q^{ik}}\sum_{jk}Q^{ij}Q^{ik}\rho_{ijk}-\frac{1}{\beta_1}\sum_{jk}Q^{ik}Q^{ij}\log Q^{ij}\nonumber\\
&\qquad\qquad\qquad\qquad\qquad\quad+\frac{1}{\beta_2}\sum_{jk}Q^{ij}Q^{ik}\log Q^{ik}\\
&\text{subject to}\ \sum_jQ^{ij}=1, \ \sum_kQ^{ik}=1,  \ Q^{ij}\geq 0,  \ Q^{ik}\geq 0,
\end{align*} 
where $Q^{ij}$$=$$\sigma_t(u_t|x_t)$, $Q^{ik}$$=$$\tau_t(w_t|x_t)$ and $\rho_{ijk}$$=$$\rho_t(x_t,u_t,w_t)$. Note that due to constraints $\sum_jQ^{ij}$$=$$1$ and $\sum_kQ^{ik}$$=$$1$, we can replace $\frac{1}{\beta_1}\sum_{jk}Q^{ik}Q^{ij}\log Q^{ij}$ and $\frac{1}{\beta_2}\sum_{jk}Q^{ij}Q^{ik}\log Q^{ik}$ by $\frac{1}{\beta_1}\sum_{j}Q^{ij}\log Q^{ij}$ and $\frac{1}{\beta_2}\sum_{k}Q^{ik}\log Q^{ik}$, respectively. For now, we neglect the non-negativity constraints and write the Lagrangian for the above optimization problem as
\begin{align*}
L=&\sum_{jk}Q^{ij}Q^{ik}\rho_{ijk}-\frac{1}{\beta_1}\sum_{j}Q^{ij}\log Q^{ij}\\
&\qquad\qquad+\frac{1}{\beta_2}\sum_{k}Q^{ik}\log Q^{ik}+\lambda^j(\sum_jQ^{ij}-1)\\
&\qquad\qquad+\lambda^k(\sum_kQ^{ik}-1)
\end{align*}
where $\lambda^j$,$\lambda^k$ are Lagrange multipliers. Then, taking derivative with respect to $Q^{ik}$ and equating it to zero, we obtain
\begin{align*}
\frac{\partial L}{\partial Q^{ik}}=\sum_{j}Q^{ij}\rho_{ijk}+\frac{1}{\beta_2}\log Q^{ik}_{\star}+\frac{1}{\beta_2}+\lambda^k=0.
\end{align*}
Rearranging terms and using the constraint $\sum_kQ^{ik}=1$, we obtain
\begin{align*}
Q^{ik}_{\star}=\frac{\exp(-\beta_2\sum_jQ^{ij}\rho_{ijk})}{\sum_{k}\exp(-\beta_2\sum_jQ^{ij}\rho_{ijk})}
\end{align*}
which is the same as \eqref{player_2_st}. Note that the resulting strategy also satisfies the non-negativity constraint. Plugging $Q^{ik}_{\star}$ into Lagrangian $L$, we obtain the optimization problem given in \eqref{player_1_st} for which the optimal variables correspond to the optimal strategy of player 1. Similarly, the optimal value \eqref{value_game_opt} of the resulting optimization problem is the value of the game. Uniqueness of the value follows from the fact that the value of the game is the optimal value of a convex optimization problem given in \eqref{player_1_st}. $\quad \Box$

 \end{document}